# CONFIDENCE SETS FOR NONPARAMETRIC WAVELET REGRESSION

By Christopher R. Genovese[1] and Larry Wasserman[2]

*Carnegie Mellon University*

We construct nonparametric confidence sets for regression functions using wavelets that are uniform over Besov balls. We consider both thresholding and modulation estimators for the wavelet coefficients. The confidence set is obtained by showing that a pivot process, constructed from the loss function, converges uniformly to a mean zero Gaussian process. Inverting this pivot yields a confidence set for the wavelet coefficients, and from this we obtain confidence sets on functionals of the regression curve.

**1. Introduction.** Wavelet regression is an effective method for estimating inhomogeneous functions. Donoho and Johnstone (1995a, b, 1998) showed that wavelet regression estimators based on nonlinear thresholding rules converge at the optimal rate simultaneously across a range of Besov and Triebel spaces. The practical implication is that, for denoising an inhomogeneous signal, wavelet thresholding outperforms linear techniques. See, for instance, Cai (1999), Cai and Brown (1998), Efromovich (1999), Johnstone and Silverman (2002) and Ogden (1997). However, confidence sets for the wavelet estimators may not inherit the convergence rate of function estimators. Indeed, Li (1989) shows that uniform nonparametric confidence sets for regression estimators decrease in radius at a $n^{-1/4}$ rate. However, with additional assumptions, Picard and Tribouley (2000) show that it is possible to get a faster rate for pointwise intervals.

In this paper we show how to construct uniform confidence sets for wavelet regression. More precisely, we construct a confidence sphere in the $\ell^2$-norm

Received August 2002; revised April 2004.
[1]Supported in part by NSF Grant SES-98-66147.
[2]Supported in part by NIH Grants R01-CA54852-07 and MH57881, and NSF Grants DMS-98-03433 and DMS-01-04016.

for the wavelet coefficients of a regression function $f$. We use the strategy of Beran and Dümbgen (1998), originating from an idea in Stein (1981), in which one constructs a confidence set by using the loss function as an asymptotic pivot. Specifically, let $\mu_1, \mu_2, \ldots$ be the coefficients for $f$ in the orthonormal wavelet basis $\phi_1, \phi_2, \ldots$, and let $(\hat{\mu}_1, \hat{\mu}_2, \ldots)$ be corresponding estimates that depend on a (possibly vector-valued) tuning parameter $\lambda$. Let $L_n(\lambda) = \sum_i (\hat{\mu}_i(\lambda) - \mu_i)^2$ be the loss function and let $S_n(\lambda)$ be an unbiased estimate of $L_n(\lambda)$. The Beran–Dümbgen strategy has the following steps:

1. Show that the *pivot process* $B_n(\lambda) = \sqrt{n}(L_n(\lambda) - S_n(\lambda))$ converges weakly to a Gaussian process with covariance kernel $K(s, t)$.
2. Show that $B_n(\hat{\lambda}_n)$ also has a Gaussian limit, where $\hat{\lambda}_n$ minimizes $S_n(\lambda)$. This step follows from the previous step if $\hat{\lambda}_n$ is independent of the pivot process or if $B_n(\hat{\lambda}_n)$ is stochastically very close to $B_n(\lambda_n)$ for an appropriate deterministic sequence $\lambda_n$.
3. Find a consistent estimator $\hat{\tau}_n^2$ of $K(\hat{\lambda}_n, \hat{\lambda}_n)$.
4. Conclude that

$$\mathcal{D}_n = \left\{\mu : \frac{L_n(\hat{\lambda}_n) - S_n(\hat{\lambda}_n)}{\hat{\tau}_n/\sqrt{n}} \leq z_\alpha\right\}$$
$$= \left\{\mu : \sum_{\ell=1}^n (\hat{\mu}_{n\ell} - \mu_\ell)^2 \leq \frac{\hat{\tau}_n z_\alpha}{\sqrt{n}} + S_n(\hat{\lambda}_n)\right\}$$

is an asymptotic $1 - \alpha$ confidence set for the coefficients, where $z_\alpha$ denotes the upper-tail $\alpha$-quantile of a standard Normal and where $\hat{\mu}_{n\ell} \equiv \hat{\mu}_\ell(\hat{\lambda}_n)$.
5. It follows that

$$\mathcal{A}_n = \left\{\sum_{\ell=1}^n \mu_\ell \phi_\ell(\cdot) : \mu \in \mathcal{D}_n\right\}$$

is an asymptotic $1 - \alpha$ confidence set for $f_n = \sum_{\ell=1}^n \mu_\ell \phi_\ell$.
6. With appropriate function-space assumptions, conclude that dilating $\mathcal{A}_n$ yields a confidence set for $f$.

The limit laws—and, thus, the confidence sets—we obtain are uniform over Besov balls. The exact form of the limit law depends on how the $\mu_i$'s are estimated. We consider universal shrinkage [Donoho and Johnstone (1995a)], modulation estimators [Beran and Dümbgen (1998)] and a restricted form of SureShrink [Donoho and Johnstone (1995b)].

Having obtained the confidence set $\mathcal{A}_n$, we immediately get confidence sets for any functional $T(f)$. Specifically, $(\inf_{f \in \mathcal{C}_n} T(f), \sup_{f \in \mathcal{C}_n} T(f))$ is an asymptotic confidence set for $T(f)$. In fact, if $\mathcal{T}$ is a set of functionals, then the collection $\{(\inf_{f \in \mathcal{C}_n} T(f), \sup_{f \in \mathcal{C}_n} T(f)) : T \in \mathcal{T}\}$ provides simultaneous intervals for all the functionals in $\mathcal{T}$. If the functionals in $\mathcal{T}$ are point-evaluators, we obtain a confidence band for $f$; see Section 8 for a discussion of confidence bands. An alternative method for constructing confidence bands is given in Picard and Tribouley (2000). At the cost of additional assumptions, the confidence set $\mathcal{A}_n$ can be expanded to a confidence set for $f$.

In Section 2 we discuss the basic framework of wavelet regression. In Section 3 we give the formulas for the confidence sets with known variance. In Section 4 we extend the results to the unknown variance case. In Section 5 we describe how to obtain confidence sets for functionals. In Section 6 we consider numerical examples. Finally, Section 7 contains technical results and Section 8 contains closing remarks.

**2. Wavelet regression.** Let $\phi$ and $\psi$ be, respectively, a father and mother wavelet that generate the following complete orthonormal set in $L^2[0,1]$:

$$\phi_{J_0,k}(x) = 2^{J_0/2}\phi(2^{J_0}x - k),$$
$$\psi_{j,k}(x) = 2^{j/2}\psi(2^j x - k),$$

for integers $j \geq J_0$ and $k$, where $J_0$ is fixed. Any function $f \in L^2[0,1]$ may be expanded as

$$(1) \qquad f(x) = \sum_{k=0}^{2^{J_0}-1} \alpha_k \phi_{J_0,k}(x) + \sum_{j=J_0}^{\infty} \sum_{k=0}^{2^j-1} \beta_{j,k} \psi_{j,k}(x),$$

where $\alpha_k = \int f \phi_{J_0,k}$ and $\beta_{j,k} = \int f \psi_{j,k}$. For fixed $j$, we call $\beta_{j,\cdot} = \{\beta_{j,k} : k = 0, \ldots, 2^j - 1\}$ the *resolution-$j$ coefficients*.

Assume that

$$Y_i = f(x_i) + \sigma \varepsilon_i, \qquad i = 1, \ldots, n,$$

where $f \in L^2[0,1]$, $x_i = i/n$ and $\varepsilon_i$ are IID standard Normals. (See Section 7 for details.) The goal is to estimate $f$ under squared error loss. We assume that $n = 2^{J_1}$ for some integer $J_1$. Let

$$(2) \qquad f_n(x) = \sum_{k=0}^{2^{J_0}-1} \alpha_k \phi_{J_0,k}(x) + \sum_{j=J_0}^{J_1} \sum_{k=0}^{2^j-1} \beta_{j,k} \psi_{j,k}(x)$$

denote the projection of $f$ onto the span of the first $n$ basis elements.

Define empirical wavelet coefficients

$$\tilde{\alpha}_k = \sum_{i=1}^n Y_i \int_{(i-1)/n}^{i/n} \phi_{j_0,k}(x)\,dx \approx \frac{1}{n}\sum_{i=1}^n \phi_{J_0,k}(x_i) Y_i \approx \alpha_k + \frac{\sigma}{\sqrt{n}} Z_k,$$

$$\tilde{\beta}_{j,k} = \sum_{i=1}^n Y_i \int_{(i-1)/n}^{i/n} \psi_{j,k}(x)\,dx \approx \frac{1}{n}\sum_{i=1}^n \psi_{j,k}(x_i) Y_i \approx \beta_{j,k} + \frac{\sigma}{\sqrt{n}} Z_{j,k},$$

where the $Z_k$s and $Z_{j,k}$s are IID standard Normals. In practice, these coefficients are computed in $O(n)$ time using the discrete wavelet transform.

We consider two types of estimation: soft thresholding and modulation. The soft-threshold estimator with threshold $\lambda \geq 0$, given by Donoho and Johnstone (1995a, 1995b), is defined by

$$\hat{\alpha}_k = \tilde{\alpha}_k, \tag{3}$$

$$\hat{\beta}_{j,k} = \operatorname{sgn}(\tilde{\beta}_{j,k})(|\tilde{\beta}_{j,k}| - \lambda)_+, \tag{4}$$

where $a_+ \equiv \max(a, 0)$.

Two common rules for choosing the threshold $\lambda$ are the universal threshold and the SureShrink threshold. To define these, let $\hat{\sigma}^2$ be an estimate of $\sigma^2$ and let $\rho_n = \sqrt{2 \log n}$. The *universal threshold* is $\lambda = \rho_n \hat{\sigma}/\sqrt{n}$. The *levelwise SureShrink* rule chooses a different threshold $\lambda_j$ for the $n_j = 2^j$ coefficients at resolution level $j$ by minimizing Stein's unbiased risk estimator (SURE) with estimated variance. This is given by

$$S_n(\lambda) = \frac{\hat{\sigma}^2}{n} 2^{J_0} + \sum_{j=J_0}^{J_1} S(\lambda_j), \tag{5}$$

where

$$S_j(\lambda_j) = \sum_{k=1}^{n_j} \left[ \frac{\hat{\sigma}^2}{n} - 2 \frac{\hat{\sigma}^2}{n} \mathbb{1}\{|\tilde{\beta}_{j,k}| \leq \lambda_j\} + \min(\tilde{\beta}_{j,k}^2, \lambda_j^2) \right], \tag{6}$$

for $J_0 \leq j \leq J_1$. The minimization is usually performed over $0 \leq \lambda_j \leq \rho_{n_j} \hat{\sigma}/\sqrt{n}$, although we shall minimize over a smaller interval for reasons that are explained in the remark after Theorem 3.2. SureShrink can also be used to select a global threshold by minimizing $S_n(\lambda)$ using the same constant $\lambda$ at every level. We call this *global SureShrink*.

The second estimator we consider is the modulation estimator given by Beran and Dümbgen (1998) and Beran (2000). Although these papers did not explicitly consider wavelet estimators, we can adapt their technique to construct estimators of the form

$$\hat{\alpha}_k = \xi_\phi \tilde{\alpha}_k, \tag{7}$$

$$\hat{\beta}_{j,k} = \xi_j \tilde{\beta}_{j,k}, \tag{8}$$

where $\xi_\phi, \xi_{J_0}, \xi_{J_0+1}, \ldots, \xi_{J_1}$ are chosen to minimize SURE, which in this case is

$$
\begin{aligned}
\tilde{S}_n(\xi) &= \sum_{k=0}^{2^{J_0}-1} \left[ \xi_\phi^2 \frac{\hat{\sigma}^2}{n} + (1-\xi_\phi)^2 \left( \tilde{\alpha}_k^2 - \frac{\hat{\sigma}^2}{n} \right) \right] \\
&\quad + \sum_{j=J_0}^{J_1} \sum_{k=0}^{2^j-1} \left[ \xi_j^2 \frac{\hat{\sigma}^2}{n} + (1-\xi_j)^2 \left( \tilde{\beta}_{j,k}^2 - \frac{\hat{\sigma}^2}{n} \right) \right] \\
&\equiv S_\phi(\xi_\phi) + \sum_{j=J_0}^{J_1} \tilde{S}_j(\xi_j).
\end{aligned}
\tag{9}
$$

Following Beran (2000), we minimize $\tilde{S}_n(\xi)$ subject to a monotonicity constraint: $1 \geq \xi_\phi \geq \xi_{J_0} \geq \xi_{J_0+1} \geq \cdots \geq \xi_{J_1}$. We call this the monotone modulator, and we let $\hat{\xi}$ denote the $\xi$'s at which the minimum is achieved.

It is natural to consider minimizing $\tilde{S}_n(\xi)$, level by level [Donoho and Johnstone (1995a, 1995b)] or in other block minimization schemes [Cai (1999)] without the monotonicity constraint. However, we find, as in Beran and Dümbgen (1998), that the loss functions for these estimators then do not admit an asymptotic distributional limit which is needed for the confidence set. It is possible to construct other modulators besides the monotone modulator that admit a limiting distribution; we will report on these elsewhere.

Having estimated the wavelet coefficients, we then estimate $f$—more precisely, $f_n$—by

$$
\hat{f}_n(x) = \sum_{k=0}^{2^{J_0}-1} \hat{\alpha}_k \phi_{J_0,k}(x) + \sum_{j=J_0}^{J_1} \sum_{k=0}^{2^j-1} \hat{\beta}_{j,k} \psi_{j,k}(x).
\tag{10}
$$

It will be convenient to consider the wavelet coefficients, true and estimated, in the form of a single vector. Let $\mu = (\mu_1, \mu_2, \ldots)$ be the sequence of true wavelet coefficients $(\alpha_0, \ldots, \alpha_{2^{J_0}-1}, \beta_{J_0,0}, \ldots, \beta_{J_0,2^{J_0}-1}, \ldots)$. The $\alpha_k$ coefficient corresponds to $\mu_\ell$, where $\ell = k+1$ and $\beta_{jk}$ corresponds to $\mu_\ell$, where $\ell = 2^j + k + 1$. Let $\phi_1, \phi_2, \ldots$ denote the corresponding basis functions. Because $f \in L^2[0,1]$, we also have that $\mu \in \ell^2$. Similarly, let $\mu^n = (\mu_1, \ldots, \mu_n)$ denote the vector of first $n$ coefficients $(\alpha_0, \ldots, \alpha_{2^{J_0}-1}, \beta_{J_0,0}, \ldots, \beta_{J_0,2^{J_0}-1}, \ldots, \beta_{J_1,2^{J_1}-1})$.

For any $c > 0$, define

$$
\ell^2(c) = \left\{ \mu \in \ell^2 : \sum_{\ell=1}^{\infty} \mu_\ell^2 \leq c^2 \right\},
$$

and let $\mathcal{B}_{p,q}^\varsigma(c)$ denote a Besov space with radius $c$. If the wavelets are $r$-regular with $r > \varsigma$, the wavelet coefficients of a function $f \in \mathcal{B}_{p,q}^\varsigma(c)$ satisfy

$\|\mu\|_{p,q}^\varsigma \leq c$, where

$$\|\mu\|_{p,q}^\varsigma = \left( \sum_{j=J_0}^{\infty} \left( 2^{j(\varsigma+(1/2)-(1/p))} \left( \sum_k |\beta_{j,k}|^p \right)^{1/p} \right)^q \right)^{1/q}. \tag{11}$$

Let

$$\gamma = \begin{cases} \varsigma, & p \geq 2, \\ \varsigma + \dfrac{1}{2} - \dfrac{1}{p}, & 1 \leq p < 2. \end{cases} \tag{12}$$

We assume that $p, q \geq 1$ and also that $\gamma > 1/2$. We also assume that the mother and father wavelets are bounded, have compact support and have derivatives with finite $L^2$ norms. We will call a space of functions $f$ satisfying these assumptions a Besov ball with $\gamma > 1/2$ and radius $c$ and the corresponding body of coefficients with $\|\mu\|_{p,q}^\varsigma \leq c$ a Besov body with $\gamma > 1/2$ and radius $c$. We use $\mathcal{B}$ to denote either, depending on context. If $\mathcal{B}$ is a coefficient body, we will denote by $\mathcal{B}^m$ for any positive integer $m$, the set of vectors $(\mu_1, \ldots, \mu_m)$ for $\mu \in \mathcal{B}$.

Our main results also extend to unions of Besov balls (and bodies). Fix $\eta, c > 0$, and define

$$\mathcal{F}_{\eta, c} = \bigcup_{p, q \geq 1} \bigcup_{\gamma \geq 1/2 + \eta} \mathcal{B}_{p,q}^{\varsigma(\gamma)}(c). \tag{13}$$

The parameter $\eta$ is an increment of smoothness required only in the non-sparse case ($p \geq 2$).

**3. Confidence sets with $\sigma$ known.** Here we give explicit formulas for the confidence set when $\sigma$ is known. The proofs are deferred until Section 7, and the $\sigma$ unknown case is treated in Section 4. It is to be understood in this section that $\sigma$ replaces $\hat\sigma$ in (5) and (9).

The confidence set is of the form

$$\mathcal{D}_n = \left\{ \mu^n : \sum_{\ell=1}^n (\mu_\ell - \hat\mu_\ell)^2 \leq s_n^2 \right\}. \tag{14}$$

The definition of the radius $s_n$ is given in Theorems 3.1, 3.2 and 3.3. In each case we will show that

$$\lim_{n \to \infty} \sup_{\mu^n \in \mathcal{B}^n} |\mathsf{P}\{\mu^n \in \mathcal{D}_n\} - (1 - \alpha)| = 0 \tag{15}$$

for a coefficient body $\mathcal{B}$. Strictly speaking, the confidence set $\mathcal{D}_n$ is for approximate wavelet coefficients, but we show in Section 7 that the approximation error can be easily accounted for. By the Parseval relation, $\mathcal{D}_n$ also yields a confidence set for $f_n$. That is,

$$\lim_{n \to \infty} \sup_{\mu^n \in \mathcal{B}^n} |\mathsf{P}\{f_n \in \mathcal{A}_n\} - (1 - \alpha)| = 0, \tag{16}$$

where

$$\mathcal{A}_n = \left\{ \sum_{j=1}^n \mu_j \phi_j : \mu^n \in \mathcal{D}_n \right\}. \tag{17}$$

Constructing the confidence set $\mathcal{A}_n$ does not require knowledge of $c$ or $\gamma$.

At the cost of making an additional assumption, namely, an upper bound on the ball size, $\mathcal{A}_n$ can be dilated slightly to produce a confidence set for $f$. Fix $\eta, c > 0$ and recall the definition of $\mathcal{F}_{\eta,c}$ from (13). Then the set

$$\mathcal{C}_n = \left\{ f \in \mathcal{F}_{\eta,c} : \inf_{g \in \mathcal{A}_n} \|f - g\|_2 \leq \frac{\delta}{\sqrt{n}} \right\}, \tag{18}$$

for $\delta > 0$, satisfies

$$\liminf_{n \to \infty} \inf_{f \in \mathcal{F}_{\eta,c}} \mathsf{P}\{f \in \mathcal{C}_n\} \geq 1 - \alpha. \tag{19}$$

The factor $\delta/\sqrt{n}$ accommodates the difference between the true and approximate wavelet coefficients. The overcoverage of (18) occurs because one never really estimates $f$, rather, any data-based procedure is inevitably estimating $f_n$.

REMARK 3.1. It is not surprising that sharp inferences are available for $f_n$ only. The difference between $f$ and $f_n$ is effectively not estimable. In the context of kernel density estimation, Neumann (1998) and Chaudhuri and Marron (2000) argue that it is sensible to confine inferences to the smoothed version of the unknown density.

REMARK 3.2. The theorems that follow state that the confidence sets have correct asymptotic coverage over a Besov space $\mathcal{B}$ with $\gamma > 1/2$. These results all hold replacing $\mathcal{B}$ by $\mathcal{F}_{\eta,c}$ for any $\eta, c > 0$. It is also worth noting that if $p < 2$, the results still hold with $\gamma = 1/2$.

THEOREM 3.1 (Universal threshold). *Suppose that $\hat{f}_n$ is the estimator based on the global threshold $\lambda = \rho_n \sigma / \sqrt{n}$. Let*

$$s_n^2 = \sigma^2 \frac{z_\alpha}{\sqrt{n/2}} + S_n(\lambda). \tag{20}$$

*Then* (15), (16) *and* (19) *hold for any Besov body $\mathcal{B}$ with $\gamma > 1/2$ and radius $c > 0$.*

We consider a restricted version of the SureShrink estimator where we minimize SURE over $\varrho \rho_n \sigma / \sqrt{n} \leq \lambda \leq \rho_n \sigma / \sqrt{n}$, where $\varrho > 1/\sqrt{2}$.

THEOREM 3.2 (Restricted SureShrink). *Let $1/\sqrt{2} < \varrho < 1$. In the global case, let $\hat{\lambda}_{J_0} = \cdots = \hat{\lambda}_{J_1} \equiv \hat{\lambda}$ be obtained by minimizing $S_n(\lambda)$ over $\varrho\rho_n\sigma/\sqrt{n} \leq \lambda \leq \rho_n\sigma/\sqrt{n}$. In the levelwise case, let $\hat{\lambda} \equiv (\hat{\lambda}_{J_0}, \ldots, \hat{\lambda}_{J_1})$ be obtained by minimizing $S_n(\lambda_{J_0}, \ldots, \lambda_{J_1})$. Let*

$$s_n^2 = \sigma^2 \frac{z_\alpha}{\sqrt{n/2}} + S_n(\hat{\lambda}). \tag{21}$$

*Then (15), (16) and (19) hold for any Besov body $\mathcal{B}$ with $\gamma > 1/2$ and radius $c > 0$.*

REMARK 3.3. We conjecture that our results hold with only the restriction that $\varrho > 0$. We hope to report on this extension in a future paper. Interestingly, the above theorem does not hold for $\varrho = 0$ because the asymptotic equicontinuity of $B_n$ fails, so some restriction on SureShrink appears to be necessary.

REMARK 3.4. The theorem also holds with a data-splitting scheme similar to that used in Nason (1996) and Picard and Tribouley (2000), where we use one half of the data to estimate the SURE-minimizing threshold and the other half to construct the confidence set. In the case $\varrho > 1/\sqrt{2}$ this is not required, but it may be needed in the more general case $\varrho > 0$.

Finally, we consider the modulation estimator.

THEOREM 3.3 (Modulators). *Let $\hat{f}_n$ be the estimate obtained from the monotone modulator. Let*

$$s_n^2 = \hat{\tau}\frac{z_\alpha}{\sqrt{n}} + \tilde{S}(\hat{a}), \tag{22}$$

*where*

$$\hat{\tau}^2 = \frac{2\sigma^4}{n}\sum_{\ell=1}^{n}(2\hat{\xi}_\ell - 1)^2 + 4\sigma^2 \sum_{\ell=1}^{n}\left(\tilde{\mu}_\ell^2 - \frac{\sigma^2}{n}\right)^2(1 - \hat{\xi}_\ell)^2, \tag{23}$$

*where $\hat{\xi}_\ell$ is the estimated shrinkage coefficient associated with $\mu_\ell$. Then (15), (16) and (19) hold for any Besov body $\mathcal{B}$ with $\gamma > 1/2$ and radius $c > 0$.*

**4. Confidence sets with $\sigma$ unknown.** Suppose now that $\sigma$ is not known. We consider two cases. The first, assumed in Beran and Dümbgen [(1998), equation 3.2], is that there exists an independent, uniformly consistent estimate of $\sigma$. For example, if there are replications at each design point, then the residuals at these points provide the required estimator $\hat{\sigma}$. More generally, letting $\mathcal{L}(\cdot)$ denote the law of a random variable, they assume the following condition:

(S1) There exists an estimate $\hat{\sigma}_n^2$, independent of the empirical wavelet coefficients, such that $\mathcal{L}(\hat{\sigma}_n^2/\sigma^2)$ depends only on $n$ and such that

$$\lim_{n \to \infty} m(\mathcal{L}(n^{1/2}(\hat{\sigma}_n^2/\sigma^2 - 1)), N(0, \mho^2)) = 0,$$

where $m(\cdot, \cdot)$ metrizes weak convergence and $\mho > 0$.

In the absence of replication (or any other independent estimate of $\sigma^2$), we estimate $\sigma^2$ by

$$\hat{\sigma}_n^2 = 2 \sum_{\ell=(n/2)+1}^{n} \tilde{\mu}_\ell^2, \tag{24}$$

which Beran (2000) calls the high-component estimator. We then need to assume that $\mu^n$ is contained in a more restrictive space. Specifically, we assume the following:

(S2) The coefficients $\mu$ of $f$ are contained in the set

$$\{\mu \in \ell^2(c) : \|\beta_{j\cdot}\|^2 \leq \zeta_j, j \geq J_2\}$$

for some $c > 0$, $J_2 > J_0$ and some sequence of positive reals $\zeta = (\zeta_1, \zeta_2, \dots)$, where $\zeta_j = O(2^{-j/2})$ and $\beta_{j\cdot}$ denotes the resolution-$j$ coefficients.

Condition (S2) holds when $f$ is in a Besov ball $\mathcal{B}$ with $\gamma > 1/2$. We note that such a condition is implicit in Beran (2000) and Beran and Dümbgen (1998) in the absence of (S1).

Beran and Dümbgen (1998) construct confidence sets with $\sigma$ unknown by including an extra term in the formula for $s_n^2$ to account for the variability in $\hat{\sigma}_n^2$. This strategy is feasible for modulators since terms involving $\hat{\sigma}_n^2$ separate nicely in the estimated loss from the rest of the data. In thresholding estimators the empirical process in Theorem 7.2 depends on $\hat{\sigma}_n$ in a complicated way, making it difficult to deal with $\hat{\sigma}$ separately. We offer two methods for this case. For the soft-thresholded wavelet estimators it turns out that a plug-in method suffices. More generally, we can use a "double confidence set" approach.

For both approaches we need the uniform consistency of $\hat{\sigma}$.

LEMMA 4.1. *For any Besov body $\mathcal{B}$ with $\gamma > 1/2$ and for every $\varepsilon > 0$,*

$$\sup_{\mu \in \mathcal{B}} \mathsf{P}\left\{ \left| \frac{\hat{\sigma}^2}{\sigma^2} - 1 \right| > \varepsilon \right\} \to 0. \tag{25}$$

The proof of this lemma is straightforward and is omitted.

In the plug-in approach we simply replace $\sigma$ by $\hat{\sigma}$ in the expressions of the last section.

THEOREM 4.1 (Plug-in confidence ball). *Theorems* 3.1 *and* 3.2 *continue to hold if* $\hat{\sigma}$ *replaces* $\sigma$. *For the modulation estimator Theorem* 3.3 *holds with* $\hat{\tau}^2$ *replaced by*

$$
\begin{aligned}
(26) \quad \hat{\tau}^2 = {} & \frac{2\hat{\sigma}^4}{n}\sum_{\ell=1}^{n}(2\hat{\xi}_\ell - 1)^2 + 2\mho\hat{\sigma}^4\left(\frac{1}{n}\sum_{\ell=1}^{n}(2\hat{\xi}_\ell - 1)\right)^2 \\
& + 4\sigma^2\sum_{\ell=1}^{n}\left(\tilde{\mu}_\ell^2 - \frac{\sigma^2}{n}\right)^2(1-\hat{\xi}_\ell)^2.
\end{aligned}
$$

In the double confidence set approach, the confidence set is the "tube" equal to the union of confidence balls obtained by treating $\sigma$ as known for every value in a confidence interval for $\sigma$. We first need a uniform confidence interval for $\sigma$. This is given in the following theorem; the proof is straightforward and is omitted.

THEOREM 4.2. *Let*

$$
(27) \quad \mathcal{Q}_n = \hat{\sigma}_n^2\left[\left(1 - \frac{\mho z_{1-\alpha/2}}{\sqrt{n}}\right)^{-1}, \left(1 - \frac{\mho z_{\alpha/2}}{\sqrt{n}}\right)^{-1}\right].
$$

*Under condition* (S1) *we have*

$$
(28) \quad \liminf_{n\to\infty} \inf_{\sigma>0} \mathsf{P}\{\sigma \in \mathcal{Q}_n\} \geq 1 - \alpha.
$$

*Under condition* (S2) *with* $\mho = 2$, *we have, for any Besov body* $\mathcal{B}$ *with* $\gamma > 1/2$,

$$
(29) \quad \liminf_{n\to\infty} \inf_{\mu\in\mathcal{B},\sigma>0} \mathsf{P}\{\sigma \in \mathcal{Q}_n\} \geq 1 - \alpha.
$$

THEOREM 4.3 (Double confidence set). *Let* $\tilde{\alpha} = 1 - \sqrt{1-\alpha}$ *if* (S1) *holds and let* $\tilde{\alpha} = \alpha/2$ *if* (S2) *holds. Let* $\mathcal{Q}_n$ *be an asymptotic* $1 - \tilde{\alpha}$ *confidence interval for* $\sigma$, *as in Theorem* 4.2. *Let*

$$
(30) \quad \mathcal{D}_n = \bigcup_{\sigma \in \mathcal{Q}_n} \mathcal{D}_{n,\sigma},
$$

*where* $\mathcal{D}_{n,\sigma}$ *is a* $1 - \tilde{\alpha}$ *confidence ball for* $\mu$ *from the previous section obtained with fixed* $\sigma$. *Then*

$$
(31) \quad \liminf_{n\to\infty} \inf_{\mu^n\in\mathcal{B}^n} \mathsf{P}\{\mu^n \in \mathcal{D}_n\} \geq 1 - \alpha.
$$

*Finally, under condition* (S1) *or* (S2), *Theorems* 3.1, 3.2 *and* 3.3 *continue to hold with* (31) *replacing* (15) *and* $\mathcal{D}_n$ *as in* (30).

**5. Confidence sets for functionals.** Let $f \mapsto f_n^\star$ be the operation that takes $f$ to the approximation defined in (44). The reader can think of $f_n^\star$ as simply the projection $f_n$ of $f$ onto the span of the first $n$ basis functions. Define $\mathcal{C}_n^\star$ to be the set of $f_n^\star$ corresponding to coefficient sequences $\mu^n \in \mathcal{D}_n$. For real-valued functionals $T$, define

$$J_n^\star(T) = \left( \inf_{f_n^\star \in \mathcal{C}_n^\star} T(f_n^\star), \sup_{f_n^\star \in \mathcal{C}_n} T(f_n^\star) \right). \tag{32}$$

We then have the following immediately from the asymptotic coverage of the confidence set.

LEMMA 5.1. *Let $\mathcal{T}$ be a set of real-valued functionals on a Besov ball $\mathcal{B}$ with $\gamma > 1/2$ and radius $c > 0$. Then*

$$\liminf_{n \to \infty} \inf_{f \in \mathcal{B}} \mathsf{P}\{T(f_n^\star) \in J_n^\star(T) \text{ for all } T \in \mathcal{T}\} \geq 1 - \alpha. \tag{33}$$

We can extend the previous result to include sets of functionals of slowly increasing resolution. Let $\mathcal{F}$ be a function class and let $\mathcal{T}_n$ be a sequence of sets of real-valued functionals on $\mathcal{F}$. Define the worst-case approximation error over $\mathcal{F}$ and $\mathcal{T}_n$ by

$$r_n(\mathcal{F}, \mathcal{T}_n) = \sup_{T \in \mathcal{T}_n} \sup_{f \in \mathcal{F}} |T(f) - T(f_n^\star)|.$$

For a sequence $w_n$, define

$$J_n(T) = \left( \inf_{f_n^\star \in \mathcal{C}_n^\star} T(f_n^\star) - w_n, \sup_{f_n^\star \in \mathcal{C}_n^\star} T(f_n^\star) + w_n \right). \tag{34}$$

THEOREM 5.1. *For a function class $\mathcal{F}$ and a sequence $\mathcal{T}_n$ of sets of real-valued functionals on $\mathcal{F}$, if $w_n \geq r_n(\mathcal{F}, \mathcal{T}_n)$,*

$$\liminf_{n \to \infty} \inf_{f \in \mathcal{F}} \mathsf{P}\{T(f) \in J_n(T) \text{ for all } T \in \mathcal{T}_n\} \geq 1 - \alpha. \tag{35}$$

PROOF. Follows from the triangle inequality and Lemma 5.1. □

REMARK 5.1. If the functionals in $\mathcal{T}_n$ are point evaluators $T(f) = f(x)$, then the confidence sets above yield confidence bands.

For a given compactly-supported wavelet basis, define the integer $\kappa$ to be the maximum number of basis functions within a single resolution level whose support contains any single point:

$$\kappa = \sup\{\#\{\psi_{jk}(x) \neq 0 : 0 \leq k < 2^j\} : 0 \leq x \leq 1, j \geq J_0\}.$$

Note also that $\|\psi_{jk}\|_1 = 2^{-j/2}\|\psi\|_1$. Both $\kappa$ and $\|\psi\|_1$ are finite for all the commonly used wavelets.

As an example, we consider local averages over intervals whose length decreases with $n$.

THEOREM 5.2. *Fix a decreasing sequence $\Delta_n > 0$ and define*

$$\mathcal{T}_n = \left\{ T : T(f) = \frac{1}{b-a} \int_a^b f\, dx, 0 \leq a < b \leq 1, |b-a| \geq \Delta_n \right\}.$$

*Fix $\eta, c > 0$ and let $\mathcal{F}_{\eta,c}$ be the union of Besov balls defined in* (13).

*If the mother and father wavelets are compactly supported with $\kappa < \infty$ and $\|\psi\|_1 < \infty$ and if $\Delta_n^{-1} = o(n^\zeta/(\log n)^{\lfloor \zeta \rfloor})$ for some $0 \leq \zeta \leq 1$, then*

(36) $$r_n(\mathcal{F}_{\eta,c}, \mathcal{T}_n) = o(n^{\zeta-1}/(\log n)^{\lfloor \zeta \rfloor}).$$

*Hence, for any sequence $w_n \geq 0$ that satisfies $w_n \to 0$ and $\liminf_{n\to\infty} w_n n^{1-\zeta} \times (\log n)^{\lfloor \zeta \rfloor} > 0$,*

(37) $$\liminf_{n\to\infty} \inf_{f \in \mathcal{F}_{\eta,c}} \mathsf{P}\{T(f) \in J_n(T) \text{ for all } T \in \mathcal{T}_n\} \geq 1 - \alpha.$$

**6. Numerical examples.** Here we study the confidence sets for the zero function $f_0(x) \equiv 0$ and for the two examples considered in Beran and Dümbgen (1998). We also compare the wavelet confidence sets to confidence sets obtained from a cosine basis as in Beran (2000).

The two functions, defined on $[0,1]$, are given by

(38) $$f_1(x) = 2(6.75)^3 x^6 (1-x)^3,$$

(39) $$f_2(x) = \begin{cases} 1.5, & \text{if } 0 \leq x < 0.3, \\ 0.5, & \text{if } 0.3 \leq x < 0.6, \\ 2.0, & \text{if } 0.6 \leq x < 0.8, \\ 0.0, & \text{otherwise.} \end{cases}$$

Tables 1 and 2 report the results of a simulation using $\alpha = 0.05$, $n = 1024$, $\sigma = 1$ and 5000 iterations (which gives a 95% confidence interval for the estimated coverage of length no more than 0.025). For comparison, the radius of the standard 95% $\chi^2$ confidence ball, which uses no smoothing, is 1.074. We used a symmlet 8 wavelet basis, and all the calculations were done using the S+Wavelets package.

**7. Technical results.** Recall that the model is

$$Y_i = f(x_i) + \sigma \varepsilon_i,$$

where the $\varepsilon_i \sim N(0,1)$ are IID and $f(x) = \sum_j \mu_j \phi_j(x)$. Let $X_j$ denote the empirical wavelet coefficients given by

$$X_j = \sum_{i=1}^n Y_i \int_{(i-1)/n}^{i/n} \phi_j(x)\, dx.$$

TABLE 1
*Coverage and average confidence ball radius, by method, in the $\sigma$-known case. Here $n = 1024$ and $\sigma = 1$*

| Method | Function | Coverage | Average radius |
|---|---|---|---|
| Universal | $f_0$ | 0.951 | 0.274 |
| | $f_1$ | 0.949 | 0.299 |
| | $f_2$ | 0.935 | 0.439 |
| SureShrink (global) | $f_0$ | 0.946 | 0.270 |
| | $f_1$ | 0.941 | 0.292 |
| | $f_2$ | 0.937 | 0.401 |
| SureShrink (levelwise) | $f_0$ | 0.944 | 0.268 |
| | $f_1$ | 0.940 | 0.289 |
| | $f_2$ | 0.927 | 0.395 |
| Modulator (wavelet) | $f_0$ | 0.941 | 0.258 |
| | $f_1$ | 0.940 | 0.269 |
| | $f_2$ | 0.933 | 0.329 |
| Modulator (cosine) | $f_0$ | 0.931 | 0.253 |
| | $f_1$ | 0.930 | 0.259 |
| | $f_2$ | 0.905 | 0.318 |

Then $X^n = (X_1, \ldots, X_n)$ are multivariate Normal with

$$(40) \qquad \mathsf{E} X_j = \overline{\mu}_j + O(1/n), \qquad \operatorname{Var} X_j = \frac{\sigma^2}{n} + O(1/n^2),$$

uniformly over $\mathcal{B}$ [Donoho and Johnstone (1999)], where $\overline{\mu}_{n\ell} = \int \overline{f}_n \phi_\ell$. The $X_j$'s are asymptotically independent.

That the $X_j$'s are asymptotically independent poses no problem. Using the orthogonal discrete wavelet transform to define the empirical wavelet coefficients yields $\tilde{X}^n$ that are exactly independent. Donoho and Johnstone (1999) show that the means and variances of $\tilde{X}^n$ and $X^n$ are close. From this, it follows that the Kullback–Leibler distance—and, hence, the total variation distance—between the law of $\sqrt{n}(X^n - \overline{\mu}^n)$ and a $N_n(0, \sigma^2 I)$ tends to 0 uniformly, where $\overline{\mu}^n = (\overline{\mu}_1, \ldots, \overline{\mu}_n)$. In what follows, we may thus assume the $X_j$ are independent Normal$(\overline{\mu}_j, \sigma^2/n)$.

TABLE 2
*Coverage, by thresholding method, in the $\sigma$-unknown case using the Plug-in Confidence Ball. Again $n = 1024$ and $\sigma = 1$*

| Function | Universal | Sure GL | Sure LW | WaveMod | CosMod |
|---|---|---|---|---|---|
| $f_0$ | 0.961 | 0.955 | 0.954 | 0.955 | 0.999 |
| $f_1$ | 0.963 | 0.955 | 0.953 | 0.961 | 0.999 |
| $f_2$ | 0.938 | 0.940 | 0.929 | 0.951 | 0.997 |

It will be helpful to introduce some notation before proceeding with the ensuing sections. Let $\sigma_n^2 = \sigma^2/n$ and define $r_n = \rho_n \sigma/\sqrt{n}$, where $\rho_n = \sqrt{2\log n}$. Also define $\nu_{ni} = -\sqrt{n}\mu_i/\sigma$, and let $a_{ni} = \nu_{ni} - u\rho_n$ and $b_{ni} = \nu_{ni} + u\rho_n$. Note that $\sqrt{n}X_i/\sigma = \varepsilon_i - \nu_{ni}$. Define

$$I_{ni}(u) = \mathbb{1}\{|X_i| < ur_n\} = \mathbb{1}\{\nu_{ni} - u\rho_n < \varepsilon_i < \nu_{ni} + u\rho_n\} = \mathbb{1}\{a_{ni} < \varepsilon_i < b_{ni}\},$$

$$I_{ni}^+(u) = \mathbb{1}\{X_i > ur_n\} = \mathbb{1}\{\varepsilon_i > \nu_{ni} + u\rho_n\} = \mathbb{1}\{\varepsilon > b_{ni}\},$$

$$I_{ni}^-(u) = \mathbb{1}\{X_i < -ur_n\} = \mathbb{1}\{\varepsilon_i < \nu_{ni} - u\rho_n\} = \mathbb{1}\{\varepsilon < a_{ni}\},$$

$$J_{ni}(s,t) = \mathbb{1}\{sr_n < X_i < tr_n\} = \mathbb{1}\{\nu_{ni} + s\rho_n < \varepsilon_i < \nu_{ni} + t\rho_n\}.$$

For $0 \leq u \leq 1$ and $1 \leq i \leq n$, define

$$\begin{aligned}Z_{ni}(u) &= \sqrt{n}[(X_i - ur_n)\mathbb{1}\{X_i > ur_n\} + (X_i + ur_n)\mathbb{1}\{X_i < -ur_n\} - \mu_i]^2 \\ &\quad - \sqrt{n}[\sigma_n^2 - 2\sigma_n^2\mathbb{1}\{X_i^2 \leq u^2 r_n^2\} + \min(X_i^2, u^2 r_n^2)] \\ &= \frac{\sigma^2}{\sqrt{n}}[(\varepsilon_i^2 - 1)(1 - 2I_{ni}(u)) \\ &\quad + 2\nu_{ni}\varepsilon_i I_{ni}(u) - 2u\rho_n\varepsilon_i(I_{ni}^+(u) - I_{ni}^-(u))].\end{aligned}$$

(41)

Each $Z_{ni}$ represents the contribution of the $i$th observation to the pivot process and satisfies $\mathsf{E}Z_{ni}(u) = 0$ for every $0 \leq u \leq 1$. We also have that

(42) $$\begin{aligned}Z_{ni}^2(u) &= \frac{\sigma^4}{n}[(\varepsilon_i^2 - 1)^2 + 4\nu_{ni}^2\varepsilon_i^2 I_{ni}(u) + 4u^2\rho_n^2\varepsilon_i^2(1 - I_{ni}(u)) \\ &\quad - 4\nu_{ni}\varepsilon_i(\varepsilon_i^2 - 1)I_{ni}(u) - 4u\rho_n\varepsilon_i(\varepsilon_i^2 - 1)(I_{ni}^+(u) - I_{ni}^-(u))].\end{aligned}$$

The relevance of these definitions will become clear subsequently. Throughout this section $C'$ denotes a generic positive constant not depending on $n$, $\mu$ or $\varepsilon$, that may change from expression to expression.

7.1. *Absorbing approximation and projection errors.* As noted in the statements of Theorems 3.1, 3.2 and 3.3, the confidence set $\mathcal{C}_n$ for $\overline{\mu}^n$ induces a confidence set for $f$ uniformly over Besov spaces. In this section we make this precise.

Define

(43) $$\begin{aligned}\overline{f}_n(x) &= n\sum_{i=1}^n \mathbb{1}_{[(i-1)/n,i/n]}(x) \int_{(i-1)/n}^{i/n} f(t)\,dt \\ &= \sum_{\ell=1}^\infty \overline{\mu}_n \phi_j(x)\end{aligned}$$

and its projection

(44) $$f_n^\star(x) = \sum_{\ell=1}^n \overline{\mu}_n \phi_j(x).$$

THEOREM 7.1. *Fix $c, \eta > 0$. Let $\bar{\mathcal{F}}_{\eta,c}$ be the body corresponding to $\mathcal{F}_{\eta,c}$. Let $\mathcal{D}_n$ be defined by* (14) *and suppose that*

$$\liminf_{n \to \infty} \inf_{\bar{\mu}^n \in \bar{\mathcal{F}}_{\eta,c}} \mathsf{P}\{\bar{\mu}^n \in \mathcal{D}_n\} \geq 1 - \alpha.$$

*Let $\mathcal{C}_n$ be defined as in* (18). *Then*

(45) $$\liminf_{n \to \infty} \inf_{f \in \mathcal{F}_{\eta,c}} \mathsf{P}\{f \in \mathcal{C}_n\} \geq 1 - \alpha.$$

PROOF. From the results in Brown and Zhao (2001), it follows that $\|f - \bar{f}_n\|_2^2$ and $\|\bar{f}_n - f_n^\star\|_2^2 = \sum_{j=n+1}^\infty \bar{\mu}_j^2$ are both bounded, uniformly over $\mathcal{B}_{p,q}^\varsigma(c)$, by $(C \log n)/n^{2\gamma}$ for some $C > 0$ not depending on $p$, $q$ or $\varsigma$. It then follows that, for any $f \in \mathcal{F}_{\eta,c}$,

$$\|f - f_n^\star\|_2^2 \leq t(\|f - \bar{f}_n\|_2 + \|\bar{f}_n - f_n^\star\|_2)^2$$
$$\leq \frac{C \log n}{n^{1+2\eta}} \equiv k_n^2.$$

Let

$$\tilde{s}_n^2 = s_n^2 + \delta_n = \frac{\hat{\tau} z_\alpha}{\sqrt{n}} + \delta_n + S_n,$$

where $\delta_n = \delta \log n / \sqrt{n}$ for any fixed, small $\delta > 0$. Let $W_n^2 = \|\hat{f}_n - f_n^\star\|_2^2$. Then

$$\|\hat{f}_n - \bar{f}_n\|_2^2 = \|\hat{f}_n - f_n^\star\|_2^2 + \|f_n^\star - \bar{f}_n\|_2^2 = W_n^2 + k_n^2$$

and

$$\|f - \hat{f}_n\|_2 \leq \|f - f_n^\star\|_2 + \|\hat{f}_n - f_n^\star\|_2 \leq W_n + k_n$$

uniformly over $\mathcal{F}_{\eta,c}$. Hence,

$$\mathsf{P}\{\|\bar{f}_n - \hat{f}_n\|_2^2 > \tilde{s}_n^2\} \leq \mathsf{P}\{W^2 > \tilde{s}_n^2 - k_n^2\}$$
$$= \mathsf{P}\{W^2 > s_n^2 + \delta_n - k_n^2\}.$$

Now, $\liminf_{n \to \infty} \delta_n - k_n^2 > 0$ and so

$$\limsup_{n \to \infty} \sup_{f \in \mathcal{F}_{\eta,c}} \mathsf{P}\{W^2 > s_n^2 + \delta_n - k_n^2\} \leq \limsup_{n \to \infty} \sup_{f \in \mathcal{F}_{\eta,c}} \mathsf{P}\{W^2 > s_n^2\} \leq \alpha.$$

To do the same for $f$ we note that

$$\|\hat{f}_n - f\|_2^2 = \|\hat{f}_n - f_n^\star\|_2^2 + \|f - f_n^\star\|_2^2 + 2\langle \hat{f}_n - f_n^\star, f_n^\star - f \rangle$$
$$= \|\hat{f}_n - f_n^\star\|_2^2 + \|f - f_n^\star\|_2^2 + 2\langle \hat{f}_n - f_n^\star, f_n^\star - f_n \rangle$$
$$= \|\hat{f}_n - f_n^\star\|_2^2 + \|f - f_n^\star\|_2^2 + 2\sum_{i=1}^n (\hat{\mu}_\ell - \bar{\mu}_\ell)(\mu_\ell - \bar{\mu}_\ell)$$

$$\leq \|\hat{f}_n - f_n^\star\|_2^2 + \|f - f_n^\star\|_2^2 + 2\|\hat{f}_n - f_n^\star\|_2 \|f_n - f_n^\star\|_2$$
$$= (\|\hat{f}_n - f_n^\star\|_2 + \|f_n - f_n^\star\|_2)^2 + \|f - f_n\|_2^2$$
$$\leq (W_n + k_n)^2 + k_n^2,$$

where the last inequality follows from the results in Brown and Zhao (2001) since $\|f_n - f_n^\star\| \leq \|f - \bar{f}_n\|$. We have

$$\mathsf{P}\{\|f - \hat{f}_n\| > \tilde{s}_n^2\} \leq \mathsf{P}\{(W_n + k_n)^2 > \tilde{s}_n^2\}$$
$$= \mathsf{P}\{(W_n + k_n)^2 > s_n^2 + \delta_n\}$$
$$= \mathsf{P}\{W_n^2 + 2k_n W_n + k_n^2 > s_n^2 + \delta_n\}$$
$$\leq \mathsf{P}\{W_n^2 > s_n^2\} + \mathsf{P}\{2k_n W_n + k_n^2 > \delta_n\}.$$

The lim sup of the first term is bounded above by $\alpha$. For the second term,

$$\limsup_{n \to \infty} \sup_{f \in \mathcal{F}_{\eta,c}} \mathsf{P}\{2k_n W_n + k_n^2 > \delta_n\}$$
$$= \limsup_{n \to \infty} \sup_{f \in \mathcal{F}_{\eta,c}} \mathsf{P}\left\{W_n > \frac{\delta_n - k_n^2}{2k_n}\right\}$$
$$= \limsup_{n \to \infty} \sup_{f \in \mathcal{F}_{\eta,c}} \mathsf{P}\left\{W_n > \frac{\delta n^\eta}{2\sqrt{2C}} - \frac{\sqrt{C \log n}}{2n^{(1/2)+\eta}}\right\}$$
$$\to 0.$$

Hence, $\limsup_{n \to \infty} \mathsf{P}\{\|f - \hat{f}_n\| > \tilde{s}_n^2\} \leq \alpha$. □

7.2. *The pivot process.* In the rest of this section, for convenience, we will denote $\bar{\mu}_j$ simply by $\mu_j$. We now focus on the confidence set $\mathcal{D}_n$ for $\mu^n$ defined by

$$\mathcal{D}_n = \left\{\mu^n : \sum_{i=1}^n (\hat{\mu}_i - \mu_i)^2 \leq s_n^2\right\}.$$

Our main task in showing that $\mathcal{D}_n$ has correct asymptotic coverage is to show that the pivot process has a tight Gaussian limit. See van der Vaart and Wellner (1996) for the definition of a tight, Gaussian limit.

For $i = 1, \ldots, n$, let $j(i)$ denote the resolution level to which index $i$ belongs, and for $j = J_0, \ldots, J_1$, let $\mathcal{I}_j$ denote the set of indices at resolution level $j$, which contains $n_j = 2^j$ elements. Let $t$ be a sequence of thresholds with one component per resolution level starting at $J_0$, where each $t_j$ is in the range $[\varrho \rho_n \sigma_n, \rho_n \sigma_n]$. It is convenient to write $t = u\rho_n \sigma/\sqrt{n}$, where $u$ is a corresponding sequence of values in $[\varrho, 1]$. In levelwise thresholding, the $t_j$'s (and $u_j$s) are allowed to vary independently. In global thresholding, all of the $t_j$'s (and $u_j$s) are equal; in this case, we treat $t$ (and $u$) interchangeably as a sequence or scalar as convenient.

The soft threshold estimator $\hat{\mu}$ is defined by

$$\hat{\mu}_i(t) = (X_i - t_{j(i)})\mathbb{1}\{X_i > t_{j(i)}\} + (X_i + t_{j(i)})\mathbb{1}\{X_i < -t_{j(i)}\}, \qquad (46)$$

for $i = 1, \ldots, n$. The corresponding loss as a function of threshold is

$$L_n(t) = \sum_{i=1}^{n}(\hat{\mu}_i(t) - \mu_i)^2.$$

We can write Stein's unbiased risk estimate as

$$S_n(t) = \sum_{i=1}^{n}(\sigma_n^2 - 2\sigma_n^2\mathbb{1}\{X_i^2 \leq t_{j(i)}^2\} + \min(X_i^2, t_{j(i)}^2)) \qquad (47)$$

$$= \sum_{j=J_0}^{J_1}\sum_{i\in\mathcal{I}_j}(\sigma_n^2 - 2\sigma_n^2\mathbb{1}\{X_i^2 \leq t_j^2\} + \min(X_i^2, t_j^2)) \qquad (48)$$

$$\equiv \sum_{j=J_0}^{J_1} S_{nj}(t_j). \qquad (49)$$

In global thresholding, we will use the first expression. In levelwise thresholding, each $S_{nj}$ is a sum of $n_j$ independent terms, and the different $S_{nj}$'s are independent.

The SureShrink thresholds are defined by minimizing $S_n$. By independence and additivity, this is equivalent in the levelwise case to separately minimizing the $S_{nj}(t_j)$s over $t_j$. That is, recalling that $r_n = \rho_n\sigma/\sqrt{n}$,

$$\hat{u}_n = \arg\min_{\varrho \leq u \leq 1} S_n(u) \quad \text{and} \quad \hat{t}_n = u_n r_n \qquad \text{(global)}, \qquad (50)$$

$$\hat{u}_{nj} = \arg\min_{\varrho \leq u_j \leq 1} S_{nj}(u_j) \quad \text{and} \quad \hat{t}_{nj} = u_{nj} r_n \qquad \text{(levelwise)}. \qquad (51)$$

We now define

$$B_n(u) = \sqrt{n}(L_n(ur_n) - S_n(ur_n)). \qquad (52)$$

We regard $\{B_n(u) : u \in \mathcal{U}_\varrho\}$ as a stochastic process. Let $\varrho > 1/\sqrt{2}$. In the global case we take $\mathcal{U}_\varrho = [\varrho, 1]$. In the levelwise case we take $\mathcal{U} = [\varrho, 1]^\infty$, the set of sequences $(u_1, \ldots, u_k, 1, 1, \ldots)$ for any positive integer $k$ and any $\varrho \leq u_j \leq 1$. This process has mean zero because $S_n$ is an unbiased estimate of risk. The process $B_n$ can be written as

$$B_n(u) = \sum_{i=1}^{n} Z_{ni}(u_{j(i)}), \qquad (53)$$

where $Z_{ni}$ is defined in (41). For levelwise thresholding, $B_n(u)$ is also additive in the threshold components:

$$B_n(u) = \sum_{j=J_0}^{J_1} B_{nj}(u_j) = \sum_{j=J_0}^{J_1} \sum_{i \in \mathcal{I}_j} Z_{ni}(u_j). \tag{54}$$

Each $B_{nj}$ is of the same basic form as the sum of $n_j$ independent terms.

LEMMA 7.1. *Let $\mathcal{B}$ be a Besov body with $\gamma > 1/2$ and radius $c > 0$. The process $B_n(u)$ is asymptotically equicontinuous on $\mathcal{U}_\varrho$ uniformly over $\mu \in \mathcal{B}$ for any $\varrho > 1/\sqrt{2}$ with both global and levelwise thresholding. In fact, it is uniformly asymptotically constant in the sense that, for all $\delta > 0$,*

$$\limsup_{n \to \infty} \sup_{\mu \in \mathcal{B}} \mathsf{P}^* \left\{ \sup_{u,v \in \mathcal{U}_\varrho} |B_n(u) - B_n(v)| > \delta \right\} = 0. \tag{55}$$

PROOF. As above, let $a_{ni} = \nu_{ni} - u\rho_n$ and $b_{ni} = \nu_{ni} + u\rho_n$. From (41) we have, for $0 \leq u < v \leq 1$,

$$\begin{aligned}
&\frac{\sqrt{n}}{2\sigma^2}(Z_{ni}(u) - Z_{ni}(v)) \\
&= (\varepsilon_i^2 - 1)(I_{ni}(v) - I_{ni}(u)) - \nu_{ni}\varepsilon_i(I_{ni}(v) - I_{ni}(u)) \\
&\quad - u\rho_n\varepsilon_i(I_{ni}^+(u) - I_{ni}^-(u)) + v\rho_n\varepsilon_i(I_{ni}^+(v) - I_{ni}^-(v)) \\
&= (\varepsilon_i^2 - 1)\mathbb{1}\{u\rho_n \leq |\varepsilon_i - \nu_{ni}| < v\rho_n\} - \nu_{ni}\varepsilon_i\mathbb{1}\{u\rho_n \leq |\varepsilon_i - \nu_{ni}| < v\rho_n\} \\
&\quad - u\rho_n\varepsilon_i\mathbb{1}\{u\rho_n \leq \varepsilon_i - \nu_{ni} < v\rho_n\} + u\rho_n\varepsilon_i\mathbb{1}\{-v\rho_n \leq \varepsilon_i - \nu_{ni} < -u\rho_n\} \\
&\quad + (v-u)\rho_n\varepsilon_i(I_{ni}^+(v) - I_{ni}^-(v)) \\
&= (\varepsilon_i^2 - 1)\mathbb{1}\{u\rho_n \leq |\varepsilon_i - \nu_{ni}| < u\rho_n + (v-u)\rho_n\} \\
&\quad - b_{ni}\varepsilon_i\mathbb{1}\{b_{ni} < \varepsilon_i \leq b_{ni} + (v-u)\rho_n\} \\
&\quad - a_{ni}\varepsilon_i\mathbb{1}\{a_{ni} - (v-u)\rho_n \leq \varepsilon_i < a_{ni}\} \\
&\quad + (v-u)\rho_n\varepsilon_i[\mathbb{1}\{\varepsilon_i > b_{ni} + (v-u)\rho_n\} - \mathbb{1}\{\varepsilon_i < a_{ni} - (v-u)\rho_n\}].
\end{aligned}$$
(56)

From (56) we have that

$$\begin{aligned}
\frac{\sqrt{n}}{2\sigma^2}|Z_{ni}(u) - Z_{ni}(v)| &\leq (\varepsilon_i^2 + (|\nu_{ni}| + u\rho_n)|\varepsilon_i| + 1)\mathbb{1}\{u\rho_n \leq |\varepsilon_i - \nu_{ni}| \leq v\rho_n\} \\
&\quad + |v-u|\rho_n|\varepsilon_i|\mathbb{1}\{|\varepsilon_i - \nu_{ni}| > v\rho_n\} \\
&\leq (\varepsilon_i^2 + |\nu_{ni}||\varepsilon_i| + 1)\mathbb{1}\{u\rho_n \leq |\varepsilon_i - \nu_{ni}| \leq v\rho_n\} \\
&\quad + \rho_n|\varepsilon_i|\mathbb{1}\{|\varepsilon_i - \nu_{ni}| \geq u\rho_n\} \\
&\leq (\varepsilon_i^2 + |\nu_{ni}||\varepsilon_i| + 1)\mathbb{1}\{\varrho\rho_n \leq |\varepsilon_i - \nu_{ni}| \leq \rho_n\} \\
&\quad + \rho_n|\varepsilon_i|\mathbb{1}\{|\varepsilon_i - \nu_{ni}| \geq \varrho\rho_n\} \\
&\equiv \Delta_{ni}.
\end{aligned}$$
(57)

Let

$$\mathcal{A}_{n0} = \{1 \leq i \leq n : |\nu_{ni}| \leq 1\}, \qquad \mathcal{A}_{n1} = \{1 \leq i \leq n : 1 < |\nu_{ni}| \leq 2\rho_n\},$$

and

$$\mathcal{A}_{n2} = \{1 \leq i \leq n : |\nu_{ni}| > 2\rho_n\}.$$

Let $\mathcal{A} = \mathcal{A}_{n,1} \cup \mathcal{A}_{n,2}$, the set of $i$ such that $|\nu_{ni}| \geq 1$. Let $n_0$ be the cardinality of $\mathcal{A}$. Let $\beta = 2\gamma$ and note that $\beta > 1$ since $\gamma > 1/2$. The Besov condition implies the following:

$$\begin{aligned}
C^2 n \rho_n^2 &\geq \sum_{i=1}^{n} \nu_{ni}^2 i^\beta \geq \sum_{i \in \mathcal{A}} \nu_{ni}^2 i^\beta \\
&\geq \sum_{i \in \mathcal{A}} i^\beta \geq \sum_{i=1}^{n_0} i^\beta \\
&\geq C_2 n_0^{1+\beta},
\end{aligned} \tag{58}$$

where the last inequality holds for large enough $n_0$. It follows from (58) that

$$n_0(n) \leq C n^{1/(1+2\gamma)} \rho_n^{2/(1+2\gamma)}, \tag{59}$$

which is $o(\sqrt{n})$.

From the above, we have in the global thresholding case that

$$\begin{aligned}
&\sup_{\varrho \leq u \leq v \leq 1} |B_n(u) - B_n(v)| \\
&\leq \sup_{\varrho \leq u \leq v \leq 1} \sum_{i=1}^{n} |Z_{ni}(u) - Z_{ni}(v)| \\
&\leq \frac{2\sigma^2}{\sqrt{n}} \sum_{i=1}^{n} [(\varepsilon_i^2 + |\nu_{ni}||\varepsilon_i| + 1) \mathbb{1}\{\varrho \rho_n \leq |\varepsilon_i - \nu_{ni}| \leq \rho_n\} \\
&\qquad + \rho_n |\varepsilon_i| \mathbb{1}\{|\varepsilon_i - \nu_{ni}| \geq \varrho \rho_n\}].
\end{aligned} \tag{60}$$

We break the sum $\sum_{i=1}^{n}$ into three sums, $\sum_{i \in \mathcal{A}_{n0}} + \sum_{i \in \mathcal{A}_{n1}} + \sum_{i \in \mathcal{A}_{n2}}$, and consider these one at a time.

For the case where $|\nu_{ni}| \leq 1$, we have the following:

$$\begin{aligned}
&\frac{2\sigma^2}{\sqrt{n}} \sum_{i \in \mathcal{A}_{n0}} [(\varepsilon_i^2 + |\nu_{ni}||\varepsilon_i| + 1) \mathbb{1}\{\varrho \rho_n \leq |\varepsilon_i - \nu_{ni}| \leq \rho_n\} \\
&\qquad + \rho_n |\varepsilon_i| \mathbb{1}\{|\varepsilon_i - \nu_{ni}| \geq \varrho \rho_n\}] \\
&\leq \frac{2\sigma^2}{\sqrt{n}} \sum_{i \in \mathcal{A}_{n0}} (\varepsilon_i^2 + (1 + \rho_n)|\varepsilon_i| + 1) \mathbb{1}\{|\varepsilon_i| \geq \varrho \rho_n - 1\}.
\end{aligned}$$

Let $t_n = \varrho \rho_n - 1$. By (72) and (73), the expected value of each summand is

$$\begin{aligned}
\mathsf{E}(\varepsilon_i^2 + (1 + \rho_n)|\varepsilon_i| + 1) \mathbb{1}\{|\varepsilon_i| \geq \varrho \rho_n - 1\} \\
= 2(t_n + \rho_n + 1)\phi(t_n) + 4(1 - \Phi(t_n)) \\
= o(n^{-1/2}).
\end{aligned}$$

The entire sum thus goes to zero as well. To see the last equality, note that there exists $\delta > 0$ such that

$$\phi(t_n) = \frac{1}{\sqrt{2\pi}} \exp\left\{-\frac{1}{2}t_n^2\right\} = \frac{1}{\sqrt{2\pi e}} e^{-\varrho^2 \rho_n^2/2} e^{\varrho\rho_n}$$

$$= \frac{1}{\sqrt{2\pi e}} n^{(\sqrt{2}\varrho/\sqrt{\log n}) - \varrho^2} = o(n^{-1/2-\delta}),$$

because $\frac{\sqrt{2}\varrho}{\sqrt{\log n}} - \varrho^2 < -1/2 - \delta$ for large enough $n$, where $\delta = |\varrho^2 - 1/2|/2$.

It follows that $\rho_n \phi(t_n) = o(n^{-1/2})$, and similarly for $(1 - \Phi(t_n)) \sim \phi(t_n)/t_n$.

For the case where $1 < |\nu_{ni}| \leq 2\rho_n$, we have the following:

$$\frac{2\sigma^2}{\sqrt{n}} \sum_{i \in \mathcal{A}_{n1}} [(\varepsilon_i^2 + |\nu_{ni}||\varepsilon_i| + 1)\mathbb{1}\{\varrho\rho_n \leq |\varepsilon_i - \nu_{ni}| \leq \rho_n\}$$

$$+ \rho_n |\varepsilon_i| \mathbb{1}\{|\varepsilon_i - \nu_{ni}| \geq \varrho\rho_n\}]$$

$$\leq \frac{2\sigma^2}{\sqrt{n}} \sum_{i \in \mathcal{A}_{n1}} (\varepsilon_i^2 + 3\rho_n |\varepsilon_i| + 1)\mathbb{1}\{|\varepsilon_i - \nu_{ni}| \geq \varrho\rho_n\}.$$

The expected value of each summand is bounded by $2 + 3\rho_n$. The expected value of the entire sum is thus bounded by

$$\frac{n_0(n)}{\sqrt{n}} 2\sigma^2 (2 + 3\rho_n) \to 0,$$

because $n_0(n)\rho_n/\sqrt{n} \to 0$.

For the case where $2\rho_n < |\nu_{ni}|$, we have the following from (57):

$$\frac{2\sigma^2}{\sqrt{n}} \sum_{i \in \mathcal{A}_{n2}} [(\varepsilon_i^2 + |\nu_{ni}||\varepsilon_i| + 1)\mathbb{1}\{\varrho\rho_n \leq |\varepsilon_i - \nu_{ni}| \leq \rho_n\}$$

$$+ \rho_n |\varepsilon_i| \mathbb{1}\{|\varepsilon_i - \nu_{ni}| \geq \varrho\rho_n\}]$$

$$\leq \frac{2\sigma^2}{\sqrt{n}} \left( \sum_{i \in \mathcal{A}_{n2}} (\varepsilon_i^2 + 2\rho_n |\varepsilon_i| + 1) + \sum_{i \in \mathcal{A}_{n2}} (|\nu_{ni}| - \rho_n)|\varepsilon_i| \mathbb{1}\{|\varepsilon_i| \geq |\nu_{ni}| - \rho_n\} \right).$$

The expected value of the summands in the first term is bounded by $2 + 2\rho_n$. The expected value of the summands in the second term is bounded by $2(|\nu_{ni}| - \rho_n)\phi(|\nu_{ni}| - \rho_n)$. Hence, the expected value of the entire sum is bounded by

$$\frac{n_0(n)}{\sqrt{n}} 2\sigma^2 (2 + \rho_n + 2(|\nu_{ni}| - \rho_n)\phi(|\nu_{ni}| - \rho_n)) \to 0,$$

because $\gamma > 1/2$ implies $n_0(n)\rho_n/\sqrt{n} \to 0$.

We have shown that $\mathsf{E}\sup_{\varrho \leq u \leq v \leq 1}|B_n(u) - B_n(v)| \to 0$. The result follows for all $\delta > 0$ by Markov's inequality.

Next, consider the levelwise thresholding case. The product space $\mathcal{U}_\varrho = [\varrho, 1]^\infty$ is the set of sequences $(u_1, \ldots, u_k, 1, 1, \ldots)$ over positive integers $k$ and $\varrho \leq u_j \leq 1$. By Tychonoff's theorem, this space is compact and thus totally bounded, so $\mathcal{U}_\varrho$ is totally bounded under the product metric $d(\underline{u}, \underline{v}) = \sum_{\ell=J_0}^\infty 2^{-\ell}|u_\ell - v_\ell|$. For $\underline{u} \in \mathcal{U}^\infty$, define

$$B_n(\underline{u}) = \sum_{i=1}^n Z_{ni}(\underline{u}_{j(i)}).$$

It follows then that, for any $\underline{u}, \underline{v} \in \mathcal{U}^\infty$, $d(\underline{u}, \underline{v}) \leq 1 - \varrho$ and

$$|B_n(\underline{u}) - B_n(\underline{v})| \leq \sum_{i=1}^n |Z_{ni}(\underline{u}_{j(i)}) - Z_{ni}(\underline{v}_{j(i)})| \tag{61}$$

$$\leq \sum_{i=1}^n \sup_{u,v \in \mathcal{U}_\varrho} |Z_{ni}(u) - Z_{ni}(v)| \tag{62}$$

$$\leq \sum_{i=1}^n \Delta_{ni}, \tag{63}$$

where $\Delta_{ni}$ is the $u, v$ independent bound established above in (57). The result above shows that

$$\mathsf{E} \sup_{\underline{u}, \underline{v} \in \mathcal{U}_\varrho} |B_n(\underline{u}) - B_n(\underline{v})| \to 0. \tag{64}$$

This implies that $B_n$ is asymptotically constant (and thus equicontinuous) on $\mathcal{U}_\varrho$. □

LEMMA 7.2. *Let $\mathcal{B}$ be a Besov body with $\gamma > 1/2$ and radius $c > 0$. For any fixed $u_1, \ldots, u_k$ in either global or levelwise thresholding, the vector $(B_n(u_1), \ldots, B_n(u_k))$ converges in distribution to a mean zero Gaussian on $\mathbb{R}^k$, uniformly over $\mu \in \mathcal{B}$, in the sense that*

$$\sup_{\mu \in \mathcal{B}} m(\mathcal{L}(B_n(u_1), \ldots, B_n(u_k)), N(0, \Sigma(u_1, \ldots, u_k; \mu))) \to 0,$$

*where $m$ is any metric on $\mathbb{R}^k$ that metrizes weak convergence and where $\Sigma$ represents a limiting covariance matrix, possibly different for each $\mu$.*

PROOF. We begin by showing that the Lindeberg condition holds uniformly over $\mu \in \mathcal{B}$ and over $0 \leq u \leq 1$.

First consider global thresholding. Define $\|Z_{ni}\| = \sup_{0 \le u \le 1} |Z_{ni}(u)|$. Recall that $\mathsf{E}Z_{ni} = 0$ for all $n$ and $i$. Now by (41) and (42),

$$Z_{ni}^2(u) \le \frac{2\sigma^4}{n}[(\varepsilon_i^2 - 1)^2 + 4u^2\rho_n^2\varepsilon_i^2(1 - I_{ni}(u)) + 4\nu_{ni}^2\varepsilon_i^2 I_{ni}(u)]$$

$$\equiv \aleph_1 + \aleph_2 + \aleph_3.$$

Note that none of $\aleph_1, \aleph_2$ or $\aleph_3$ depends on $u$. Hence,

(65)
$$\begin{aligned}\|Z_{ni}\|^2 \mathbb{1}\{\|Z_{ni}\| > \eta\} \\ \le (\aleph_1 + \aleph_2 + \aleph_3)\mathbb{1}\{(\aleph_1 + \aleph_2 + \aleph_3) > \eta^2\} \\ \le \sum_{r=1}^{3}\sum_{s=1}^{3} \aleph_r J_s,\end{aligned}$$

where $J_s = \mathbb{1}\{\aleph_s > \eta^2/3\}$. We will now show that the nine terms in (65) are exponentially small in $n$, which implies that the Lindeberg condition holds.

First,

$$\mathsf{P}\left\{\aleph_1 > \frac{\eta^2}{3}\right\} = \mathsf{P}\left\{|\varepsilon_i^2 - 1| > \frac{\eta\sqrt{n}}{\sigma^2\sqrt{12}}\right\} \le 2\exp\left\{-\frac{\eta\sqrt{n}}{8\sigma^2\sqrt{12}}\right\},$$

using the fact that $\mathsf{P}\{|\chi_1^2 - 1| > t\} \le 2e^{-t(t \wedge 1)/8}$. To bound $\aleph_2$, we use Mills' ratio:

$$\mathsf{P}\left\{\aleph_2 > \frac{\eta^2}{3}\right\} \le \mathsf{P}\left\{|\varepsilon_i| > \frac{\eta}{\sigma r_n\sqrt{48}}\right\} \le 2\frac{r_n\sqrt{48}}{\eta}e^{-\eta^2/(96r_n^2)} = 2\frac{\rho\sqrt{48}}{\eta\sqrt{n}}e^{-n\eta^2/(96\rho^2)}.$$

For the third term, if $\mu_i = 0$, $\aleph_3 = 0$. If $\mu_i \ne 0$,

$$\mathsf{P}\left\{\aleph_3 > \frac{\eta^2}{3}\right\} \le \mathsf{P}\left(\{|X_i| \le r_n\} \cap \left\{\varepsilon_i^2 > \frac{\eta^2}{48\sigma^2\mu_i^2}\right\}\right) \equiv b(\mu_i).$$

An elementary calculus argument shows that $b(\mu_i) \le b(\mu_*)$, where

$$|\mu_*| = \frac{\rho_n\sigma}{2\sqrt{n}} + \frac{1}{2}\sqrt{\frac{\rho_n^2\sigma^2}{n} + \frac{4\eta}{\sqrt{48n}}}.$$

Now, for all large $n$,

$$b(\mu_*) \le \mathsf{P}\{\varepsilon > -\rho_n\sigma + \sqrt{n}|\mu_*|\}$$

$$\le \mathsf{P}\left\{\varepsilon > \frac{n^{1/4}\sqrt{\eta}}{6}\right\} \le \frac{6}{\eta\sqrt{2\pi}n^{1/4}}e^{-\eta\sqrt{n}/72}.$$

These inequalities show that, for $\eta > 0$ and for $s = 1, 2, 3$, $\mathsf{E}J_{si} \le K_1 \times \exp(-K_2 \min(\eta, \eta^2)\sqrt{n})$. Because $\sqrt{\mathsf{E}\aleph_{1i}^2} \le K_3/n$, $\sqrt{\mathsf{E}\aleph_{2i}^2} \le \bar{\rho}_n^2 K_4/n$ and $\sqrt{\mathsf{E}\aleph_{3i}^2} \le \mu_i^2 K_5$, the Cauchy–Schwarz inequality and (65) show that, for $\eta > 0$,

(66) $\mathsf{E}\sum_{i=1}^{n}\|Z_{ni}\|^2\mathbb{1}\{\|Z_{ni}\| > \eta\} \le K_6(\sigma, \bar{\rho}, c)\exp(-K_7(\sigma, \bar{\rho})\min(\eta, \eta^2)\sqrt{n}).$

Here the constants $K_j$ depend, at most, on $\sigma$. It follows that the Lindeberg condition holds uniformly by applying the Cauchy–Schwarz inequality to (65).

Write $B_n(u) \equiv B_{n;\mu}(u)$ to emphasize the dependence on $\mu$ and similarly for $Z_{ni;\mu_i}(u)$. In particular, let $B_{n;0}(u)$ denote the process with $\mu_1 = \mu_2 = \cdots = 0$. Let $\mathcal{L}_{n;\mu}(u)$ denote the law of $B_{n;\mu}(u) = \sum_{i=1}^{n} Z_{ni;\mu_i}(u)$ and let $\mathcal{N}$ denote a Normal with mean 0 and variance 2. By the triangle inequality,

$$m(\mathcal{L}_{n;\mu}(u), \mathcal{N}) \leq m(\mathcal{L}_{n;0}(u), \mathcal{N}) + m(\mathcal{L}_{n;\mu}(u), \mathcal{L}_{n;0}(u)),$$

where $m(\cdot, \cdot)$ denotes the Prohorov metric. By the uniform Lindeberg condition above, the CLT holds for $\mathcal{L}_{n;0}(u)$ and, hence, by Theorem 7.3, $m(\mathcal{L}_{n;0}(u), \mathcal{N}) \to 0$. Now we show that

$$(67) \qquad \sup_{\mu \in B} m(\mathcal{L}_{n;\mu}(u), \mathcal{L}_{n;0}(u)) \to 0.$$

Note that

$$\frac{\sqrt{n}}{2\sigma^2} |Z_{ni;\mu_i}(u) - Z_{ni;0}(u)|$$
$$= |(\varepsilon_i^2 - 1)(I_{ni;\mu_i}(u) - I_{ni;0}(u)) + \nu_{ni}\varepsilon_i I_{ni;\mu_i}(u)$$
$$- u\rho_n\varepsilon_i[(I^+_{ni;\mu_i}(u) - I^+_{ni;0}(u)) - (I^-_{ni;\mu_i}(u) - I^-_{ni;0}(u))]|.$$

This can be bounded as in the proof of Lemma 7.1 and the sum split over the same three cases $|\nu_{ni}| \leq 1$, $1 < |\nu_{ni}| \leq 2\rho_n$ and $|\nu_{ni}| > 2\rho_n$. It follows that

$$(68) \qquad \sup_{\mu \in \mathcal{B}} \mathsf{E} \sup_{\varrho \leq u \leq 1} |B_{n;\mu}(u) - B_{n;0}(u)| \leq a_n^2,$$

where $a_n \to 0$; note that $a_n$ does not depend on $u$ or $\mu$. Therefore,

$$\sup_{\mu \in \mathcal{B}} \sup_{\varrho \leq u \leq 1} \mathsf{P} |B_{n;\mu}(u) - B_{n;0}(u)| > a_n \leq \frac{a_n^2}{a_n} = a_n$$

for all large $n$. Recall that, by Strassen's theorem, if $\mathsf{P}\{|X - Y| > \varepsilon\} \leq \varepsilon$, then the marginal laws of $X$ and $Y$ are no more than $\varepsilon$ apart in Prohorov distance. Hence,

$$(69) \qquad \sup_{\mu \in \mathcal{B}} \sup_{\varrho \leq u \leq 1} m(\mathcal{L}_{n;\mu}(u), \mathcal{L}_{n;0}(u)) \leq a_n \to 0.$$

This establishes the theorem for one $u$. When $B_n(u_1, \ldots, u_k)$ is an $\mathbb{R}^k$-valued process for some fixed $k$,

$$(70) \qquad \begin{aligned} &\mathsf{E}\|B_{n;\mu}(u_1, \ldots, u_k) - B_{n;0}(u_1, \ldots, u_k)\| \\ &\qquad \leq k\mathsf{E} \sup_{\varrho \leq u \leq 1} |B_{nr;\mu}(u) - B_{nr;0}(u)|, \end{aligned}$$

so by (68) the sup of the former is bounded by $ka_n^2$. Since $k$ is fixed, the result follows. Thus, (67) holds for any finite-dimensional marginal.

The same method shows that the result also holds in the levelwise case. $\square$

THEOREM 7.2. *For any Besov body with $\gamma > 1/2$ and radius $c > 0$ and for any $1/\sqrt{2} < \varrho < 1$, there is a mean zero Gaussian process $W$ such that $B_n \rightsquigarrow W$ uniformly over $\mu \in \mathcal{B}$, in the sense that*

$$\sup_{\mu \in \mathcal{B}} m(\mathcal{L}(B_n), \mathcal{L}(W)) \to 0, \tag{71}$$

*where $m$ is any metric that metrizes weak convergence on $\ell^\infty[\varrho, 1]$.*

PROOF. The result follows from the preceeding lemmas in both the global and levelwise cases. Lemmas 7.3 and 7.2 show that the finite-dimensional distributions of the process converge to Gaussian limits. Lemma 7.1 proves asymptotic equicontinuity. It follows then that $B_n$ converges weakly to a tight Gaussian process $W$. □

7.3. *The variance and covariance of $B_n$.* Recall that $r_n = \rho_n \sigma / \sqrt{n}$, $\nu_{ni} = -\sqrt{n}\mu_i/\sigma$, $a_{ni} = \nu_{ni} - u\rho_n$ and $b_{ni} = \nu_{ni} + u\rho_n$. Also define

$$D_1(s,t) = \int_s^t \varepsilon \phi(\varepsilon) \, d\varepsilon = s\phi(s) - t\phi(t), \tag{72}$$

$$D_2(s,t) = \int_s^t \varepsilon^2 \phi(\varepsilon) \, d\varepsilon = s\phi(s) - t\phi(t) + \Phi(t) - \Phi(s), \tag{73}$$

$$D_3(s,t) = \int_s^t \varepsilon(\varepsilon^2 - 1)\phi(\varepsilon) \, d\varepsilon = (s^2 + 1)\phi(s) - (t^2 + 1)\phi(t), \tag{74}$$

$$\begin{aligned} D_4(s,t) &= \int_s^t (\varepsilon^2 - 1)^2 \phi(\varepsilon) \, d\varepsilon \\ &= 2(\Phi(t) - \Phi(s)) + s(s^2 + 1)\phi(s) - t(t^2 + 1)\phi(t). \end{aligned} \tag{75}$$

Let $K_n(u,v) = \text{Cov}(B_n(u), B_n(v))$. It follows from (42) that

$$\begin{aligned} K_n(u,u) &= \mathsf{E} Z_{ni}^2(u) \\ &= \frac{2\sigma^4}{n}[1 + 2\nu_{ni}^2 D_2(a_{ni}, b_{ni}) + 2u^2\rho_n^2(1 - D_2(a_{ni}, b_{ni})) \\ &\quad - 2\nu_{ni} D_3(a_{ni}, b_{ni}) + 2u\rho_n(D_3(-\infty, a_{ni}) - D_3(b_{ni}, \infty))] \\ &= \frac{2\sigma^4}{n}[1 + 2u^2\rho_n^2 + 2a_{ni}b_{ni} D_2(a_{ni}, b_{ni}) \\ &\quad - 2b_{ni}(a_{ni}^2 + 1)\phi(a_{ni}) + 2a_{ni}\ (b_{ni}^2 + 1)\phi(b_{ni})] \\ &= \frac{2\sigma^4}{n}[1 + 2u^2\rho_n^2 + 2a_{ni}b_{ni}(\Phi(b_{ni}) - \Phi(a_{ni})) \\ &\quad + 2b_{ni}a_{ni}^2\phi(a_{ni}) - 2a_{ni}b_{ni}^2\phi(b_{ni}) \\ &\quad - 2b_{ni}(a_{ni}^2 + 1)\phi(a_{ni}) + 2a_{ni}(b_{ni}^2 + 1)\phi(b_{ni})] \end{aligned}$$

$$= \frac{2\sigma^4}{n}[1 + \underbrace{2u^2\rho_n^2 + 2a_{ni}b_{ni}(\Phi(b_{ni}) - \Phi(a_{ni}))}_{I}$$

$$+ \overbrace{2a_{ni}\phi(b_{ni}) - 2b_{ni}\phi(a_{ni})}^{II}]$$

$$\equiv \frac{2\sigma^4}{n}[1 + I + II].$$

THEOREM 7.3. *Let $\mathcal{B}$ be a Besov ball with $\gamma > 1/2$ and radius $c > 0$. Then,*

$$\lim_{n\to\infty} \sup_{\mu\in\mathcal{B}} \left|\sum_{i=1}^n \mathsf{E} Z_{ni}^2(u) - 2\sigma^4\right| = 0.$$

PROOF. Apply Lemma 7.4 to the sum of terms $I$ and $II$. This is of the form $\frac{1}{n}\sum_{i=1}^n g_n(\nu_{ni})$, where

$$g_n(x) = 2u^2\rho_n^2 + 2(x^2 - u^2\rho_n^2)(\Phi(x + u\rho_n) - \Phi(x - u\rho_n))$$
$$+ 2(x - u\rho_n)\phi(x + u\rho_n) - 2(x + u\rho_n)\phi(x - u\rho_n).$$

We have $g_n(0) \to 0$ because $|g_n(0)| \leq 6\rho_n n^{-\varrho^2}$, and, hence, $n > 288/\varepsilon$ implies that $|g_n(0)| < \varepsilon$.

Now, if $|x| > 2\rho_n$, then by Mills' inequality $|g_n(x)| \leq C\rho_n^2$. If $|x| \leq 2\rho_n$, the same holds because each term is of order $\rho_n^2$. Hence, $\|g_n\|_\infty = O(\log n)$.

For $x$ in a neighborhood of zero,

$$|g_n(x) - g_n(0)| \leq |g'_n(\xi)||x| \qquad \text{for some } |\xi| \leq |x|$$
$$\leq \sup_{|\xi|\leq|x|} |g'_n(\xi)||x|.$$

Hence,

$$\sup_n |g_n(x) - g_n(0)| \leq |x| \sup_n \sup_{|\xi|\leq|x|} |g'_n(\xi)|.$$

By direct calculation, for $\varepsilon > 0$ and $\delta = \min(\varepsilon, 1/8)$, $\sup_{|\xi|\leq|x|} |g'_n(\xi)| \leq 1$, so $|x| \leq \delta$ implies $\sup_n |g_n(x) - g_n(0)| \leq \varepsilon$. Thus, $(g_n)$ is an equicontinuous family of functions at 0.

By Lemma 7.4, the result follows. □

LEMMA 7.3. *Let $\mathcal{B}$ be a Besov body with $\gamma > 1/2$ and radius $c > 0$. Then the function $K_n(u,v) = \mathrm{Cov}(B_n(u), B_n(v))$ converges to a well-defined limit uniformly over $\mu \in \mathcal{B}$:*

$$\lim_{n\to\infty} \sup_{\mu\in\mathcal{B}} |K_n(u,v) - 2\sigma^4| = 0.$$

PROOF. Theorem 7.3 proves the result for $u = v$. Let $0 \le u < v \le 1$. Then by (41),

$$Z_{ni}(u)Z_{ni}(v) = \frac{\sigma^4}{n}(\varepsilon^2 - 1)^2(1 - 2I_{ni}(v) + 2I_{ni}(u))$$

$$+ 2\frac{\sigma^3}{\sqrt{n}}\varepsilon(\varepsilon^2 - 1)\{vr_n I_{ni}^-(v) + ur_n I_{ni}^-(u) - \mu_i$$

$$- vr_n I_{ni}^+(v) - ur_n I_{ni}^+(u) + 3\mu I_{ni}(u)$$

$$+ 2ur_n J_{ni}(u,v) - 2ur_n J_{ni}(-v,-u)\}$$

$$+ 2\sigma^2\varepsilon^2\{2\mu^2 I_{ni}(u) + 2ur_n(\mu + r_n v)J_{ni}(u,v)$$

$$+ 2ur_n(vr_n - \mu v)J_{ni}(-v,-u)\}.$$

Let $\tilde{a}_{ni} = \nu_{ni} - v\rho$ and $\tilde{b}_{ni} = \nu_{ni} + v\rho$. We then have

$$K_n(u,v) = \mathsf{E}(Z_{ni}(u)Z_{ni}(v))$$

$$= \frac{2\sigma^4}{n}[1 - D_4(\tilde{a}_{ni}, \tilde{b}_{ni}) + D_4(a_{ni}, b_{ni})$$

$$- v\rho_n D_3(-\infty, \tilde{a}_{ni}) + u\rho_n D_3(-\infty, a_{ni}) + \nu_{ni}$$

$$- v\rho_n D_3(\tilde{b}_{ni}, \infty) - u\rho_n D_3(b_{ni}, \infty) - 3\nu_{ni} D_3(a_{ni}, b_{ni})$$

$$- 2u\rho_n D_3(b_{ni}, \tilde{b}_{ni}) - 2u\rho_n D_3(\tilde{a}_{ni}, a_{ni})$$

$$+ 2\nu_{ni}^2 D_2(a_{ni}, b_{ni}) + 2uv\rho_n(\rho_n - \nu_{ni})D_2(b_{ni}, \tilde{b}_{ni})$$

$$+ 2uv\rho_n(\rho_n + \nu_{ni})D_2(\tilde{a}_{ni}, a_{ni})].$$

The proof that this converges is essentially the same as the proof of Theorem 7.3.

□

LEMMA 7.4. *Let $\mathcal{B}$ be a Besov ball with $\gamma > 1/2$. Let $g_n$ be a sequence of functions equicontinuous at 0, with $\|g_n\|_\infty = O((\log n)^\alpha)$ for some $\alpha > 0$, and satisfying $g_n(0) \to a \in \mathbb{R}$. Then*

$$\lim_{n \to \infty} \sup_{\mu \in \mathcal{B}} \frac{1}{n} \sum_{i=1}^n g_n(\mu_i \sqrt{n}) = a.$$

PROOF. Without loss of generality, assume that $a = 0$. Let $M_n = \|g_n\|_\infty$. Fix $\varepsilon > 0$. By equicontinuity, there exists $\delta > 0$ such that $|x| < \delta$ implies $|g_n(x) - g_n(0)| < \varepsilon/4$ for all $n$. By assumption, there exists an $N$ such that $|g_n(0)| < \varepsilon/4$ for $n \ge N$. Since $\mathcal{B}$ is by assumption a Besov ball, there is a constant $C$ such that, for all $n$, $\sum_{i=1}^n \mu_i^2 i^{2\gamma} \le C^2 \log n$, for all $\mu \in \mathcal{B}$. See Cai [(1999), pages 919 and 920] for inequalities that imply this.

Let $\nu_{ni} = \mu_i \sqrt{n}$. The condition on $\mu$ implies for all $n$ that

$$\sum_{i=1}^{n} \nu_{ni}^2 i^{2\gamma} \leq C^2 n \log n.$$

Let the set of such $\nu_n$s be denoted by $\tilde{\mathcal{B}}_n$. We thus have

$$\sup_{\mu \in \mathcal{B}} \left| \frac{1}{n} \sum_{i=1}^{n} g(\mu_i \sqrt{n}) \right| \leq \sup_{\nu_n \in \tilde{\mathcal{B}}_n} \frac{1}{n} \sum_{i=1}^{n} |g(\nu_{ni})|.$$

Let

$$n_0 = \lceil C^{1/\gamma} n^{1/2\gamma} (\log n)^{1/2\gamma} / \delta^{1/\gamma} \rceil.$$

This is less than $n$ and bigger than $N$ for large $n$. Then for $i \geq n_0$ and $n \geq N$, $|\nu_{ni}| \leq \delta$ and $|g_n(\nu_{ni})| \leq \varepsilon/2$. We have

$$\frac{1}{n} \sum_{i=1}^{n} |g_n(\nu_{ni})| \leq \frac{n_0}{n} M_n + \frac{\varepsilon}{2} \frac{n - n_0 + 1}{n}$$

$$\leq n^{-(1-1/2\gamma)} (\log n)^{1/2\gamma} \frac{C^{1/\gamma} M_n}{\delta^2} + \frac{\varepsilon}{2}.$$

Thus, as soon as

$$n (\log n)^{-1/(2\gamma-1)} \geq \left( \frac{C^{1/\gamma}}{\delta^{1/\gamma}} \max\left(1, \frac{2 M_n}{\varepsilon}\right) \right)^{2\gamma/(2\gamma-1)},$$

we have

$$\sup_{\nu_n \in \tilde{\mathcal{B}}_n} \frac{1}{n} \sum_{i=1}^{n} |g_n(\nu_{ni})| \leq \varepsilon,$$

which proves the lemma. □

### 7.4. *Proofs of main theorems.*

PROOF OF THEOREM 3.1. This follows from Theorems 7.2 and 7.3. The last statement follows from Theorem 7.1. □

PROOF OF THEOREM 3.2. This follows from Theorems 7.2 and 7.3 and the fact that $B(\hat{u}) = B(1) + o_P(1)$ uniformly in $\varrho \leq \hat{u} \leq 1$, and $\mu \in \mathcal{B}$. The last statement follows from Theorem 7.1. □

PROOF OF THEOREM 3.3. This follows from Theorem 3.2 in Beran and Dümbgen (1998) and Theorem 7.1. □

PROOF OF THEOREM 4.1. Let $\hat{m} = \hat{\sigma}/\sigma$. The pivot process with $\hat{\sigma}$ "plugged in" is

$$\hat{B}_n(u) = \sqrt{n} \sum_{i=1}^n [(\mu_i - \hat{\mu}_i(ur_n\hat{m}))^2$$
$$+ [\hat{m}\sigma_n^2 - 2\hat{m}\sigma_n^2 \mathbb{1}\{X_i^2 \le u^2 r_n^2 \hat{m}^2\} + \min(X_i^2, u^2 r_n^2 \hat{m}^2)]]$$
$$= B_n(u\hat{m}) + (\hat{m}^2 - 1)\frac{\sigma^2}{n}\sum_{i=1}^n (1 - 2\{|\tilde{\beta}_{j,k}| \le ur_n\hat{m}\})$$
$$= B_n(u) + o_P(1),$$

uniformly over $u \in \mathcal{U}_\varrho$ and $\mu \in \mathcal{B}$, by Lemmas 7.1 and 4.1. The result follows. □

PROOF OF THEOREM 4.3. Let $\mu_0$ and $\sigma_0$ denote the true values of $\mu$ and $\sigma$, respectively. Then under (S1) we have

$$\mathsf{P}\mu_0^n \in \mathcal{D}_n \ge \mathsf{P}\{\sigma_0 \in \mathcal{Q}_n\}\mathsf{P}\{\mu_0^n \in \mathcal{D}_n|\sigma_0 \in \mathcal{Q}_n\}$$
$$\ge \mathsf{P}\{\sigma_0 \in \mathcal{Q}_n\}\mathsf{P}\{\mu_0^n \in \mathcal{D}_{n,\sigma_0}|\sigma_0 \in \mathcal{Q}_n\}$$
$$= \mathsf{P}\{\sigma_0 \in \mathcal{Q}_n\}\mathsf{P}\{\mu_0^n \in \mathcal{D}_{n,\sigma_0}\}.$$

Hence,

$$\liminf_{n\to\infty} \inf_{\mu^n \in \mathcal{B}^n} \mathsf{P}\{f_n \in \mathcal{D}_n\} \ge (1 - \tilde{\alpha})^2 = (1 - \alpha).$$

Under (S2),

$$\mathsf{P}\{\mu_0^n \notin \mathcal{D}_n\} = \mathsf{P}\{\mu_0^n \notin \mathcal{D}_n, \sigma_0 \notin \mathcal{Q}_n\} + \mathsf{P}\{\mu_0^n \notin \mathcal{D}_n, \sigma_0 \in \mathcal{Q}_n\}$$
$$\le \mathsf{P}\{\sigma_0 \notin \mathcal{Q}_n\} + \mathsf{P}\{\mu_0^n \notin \mathcal{D}_{n,\sigma_0}, \sigma_0 \in \mathcal{Q}_n\}$$
$$\le \mathsf{P}\{\sigma_0 \notin \mathcal{Q}_n\} + \mathsf{P}\{\mu_0^n \notin \mathcal{D}_{n,\sigma_0}\}.$$

Thus,

$$\liminf_{n\to\infty} \inf_{\mu^n \in \mathcal{B}^n} \mathsf{P}\{f_n \in \mathcal{C}_n\} \ge (1 - \tilde{\alpha} - \tilde{\alpha}) = 1 - \alpha.$$

This completes the proof.

For the final claim, note that the uniform consistency of $\hat{\sigma}$ and the asymptotic constancy of $B_n$ (Lemma 7.1) imply that $B(\hat{u}) = B(1) + o_P(1)$, uniformly in $\varrho \le \hat{u} \le 1$ and $\mu \in \mathcal{B}$. The theorem follows from Theorems 3.1, 3.2 and 3.3 and 4.3. □

PROOF OF THEOREM 5.2. For any $f \in \mathcal{F}_{\eta,c}$, we have that

(76)        $|T(f) - T(f_n^\star)| \le |T(f) - T(f_n)| + |T(f_n) - T(f_n^\star)|.$

NONPARAMETRIC WAVELET REGRESSION 29

Since $\int_a^b \psi_{jk} = 0$ whenever the support of $\psi_{jk}$ is contained in $[a,b]$, the first term is bounded by (with $C'$ denoting a possibly different constant in each expression)

$$|T(f) - T(f_n)| \leq \sum_{j=J_1+1}^{\infty} \sum_{k=0}^{2^j-1} |\beta_{jk}| \frac{1}{b-a} \left| \int_a^b \psi_{jk}(x)\,dx \right|$$

$$\leq \frac{\|\psi\|_1 \kappa C'}{\Delta_n} \sum_{j=J_1+1}^{\infty} \max |\beta_{j\cdot}| 2^{-j/2}$$

$$\leq \frac{C'}{\Delta_n} \sum_{j=J_1+1}^{\infty} \|\beta_{j\cdot}\|_2 2^{-j/2}$$

$$\leq \frac{C'}{\Delta_n} \sum_{j=J_1+1}^{\infty} 2^{-j}$$

$$= \frac{C'}{\Delta_n} 2^{-J_1}$$

$$= \frac{C'}{n\Delta_n}$$

$$= o(n^{\zeta-1}/(\log n)^{\lfloor \zeta \rfloor}).$$

For a given $0 \leq a < b \leq 1$, let $q = \sup\{1 \leq m \leq n : (m-1)/n \leq a\}$ and $r = \inf\{1 \leq m \leq n : b \leq m/n\}$. The second term in (76) is bounded by

$$|T(f_n) - T(f_n^\star)| \leq \frac{1}{b-a} \sum_{\ell=1}^{n} |\mu_\ell| \left| \int_a^b (\phi_\ell - \bar{\phi}_\ell) \right|$$

$$= \frac{1}{b-a} \sum_{\ell=1}^{n} |\mu_\ell| \left| \int_a^b \phi_\ell - \int_{(q-1)/n}^{r/n} \phi_\ell \right|$$

$$\leq \frac{1}{b-a} \sum_{\ell=1}^{n} |\mu_\ell| \left( \int_{(q-1)/n}^{a} |\phi_\ell| + \int_b^{r/n} |\phi_\ell| \right)$$

$$\leq \frac{1}{b-a} \left[ \sum_{k=0}^{2^{J_0}-1} |\alpha_k| \frac{C_0}{n} + \frac{4\kappa C_1}{n} \sum_{j=J_0}^{J_1} \max |\beta_{j\cdot}| 2^{j/2} \right]$$

$$\leq \frac{1}{b-a} \left[ \sum_{k=0}^{2^{J_0}-1} |\alpha_k| \frac{C_0}{n} + \frac{4\kappa C_1}{n} \sum_{j=J_0}^{J_1} \|\beta_{j\cdot}\|_2 2^{j/2} \right]$$

$$\leq \frac{1}{b-a} \left[ \sum_{k=0}^{2^{J_0}-1} |\alpha_k| \frac{C_0}{n} + \frac{4\kappa C_1 c}{n} (J_1 - J_0) \right]$$

$$= \frac{1}{b-a} \left[ \sum_{k=0}^{2^{J_0}-1} |\alpha_k| \frac{C_0}{n} + \frac{C' \log n}{n} \right]$$

$$\leq \frac{1}{\Delta_n} \frac{C''(1 + \log n)}{n}$$

$$= o(n^{\zeta-1}/(\log n)^{\lfloor \zeta \rfloor}).$$

It follows that $r_n(\mathcal{F}_{\eta,c}, \mathcal{T}_n) = o(n^{\zeta-1}/(\log n)^{\lfloor \zeta \rfloor})$. The result follows by Theorem 5.1. $\square$

**8. Discussion.** The expected radius of the confidence ball can be shown to be of order $n^{-1/4}$. This is not surprising since the minimax estimation rate for a Besov space is $n^{-\gamma/(2\gamma+1)}$, which approaches $n^{-1/4}$ as $\gamma$ approaches $1/2$. Moreover, Li (1989) showed that for nonparametric regression without smoothness constraints, confidence spheres for nonparametric regression cannot shrink faster than $n^{-1/4}$. Indeed, the presence of the term $\tau/\sqrt{n}$ in the squared radius of our confidence balls implies that rate cannot be faster than $n^{-1/4}$. This is consistent with the results in Low (1997) and Cai and Low (2003) that suggest confidence sets cannot be rate adaptive. Thus, while we have not shown that our confidence set $\mathcal{D}_n$ is rate optimal, we doubt that the rate can be improved. One consequence of the slow rate of the confidence set is that the arguments that favor threshold estimators over modulators no longer apply.

We have chosen to emphasize confidence balls and simultaneous confidence sets for functionals. A more traditional approach is to construct an interval of the form $\hat{f}(x) \pm w_n$, where $\hat{f}(x)$ is an estimate of $f(x)$ and $w_n$ is an appropriate sequence of constants. This corresponds to taking $T(f) = f(x)$, the evaluation functional, in Theorem 5.1. There is a rich literature on this subject; a recent example in the wavelet framework is Picard and Tribouley (2000). Such confidence intervals are pointwise in two senses. First, they focus on the regression function at a particular point $x$, although they can be extended into a confidence band. Second, the validity of the asymptotic coverage usually only holds for a fixed function $f$: the absolute difference between the coverage probability and the target $1 - \alpha$ converges to zero for each fixed function, but the supremum of this difference over the function space need not converge. Moreover, in this approach one must estimate the asymptotic bias of the function estimator or eliminate the bias by undersmoothing. While acknowledging that this approach has some appeal and is certainly of great value in some cases, we prefer the confidence ball approach for several reasons. First, it avoids having to estimate and correct for the bias which is often difficult to do in practice and usually entails putting extra assumptions on the functions. Second, it produces confidence sets that are asymptotically uniform over large classes. Third, it leads directly to confidence sets for classes of functionals which we believe are quite useful in scientific applications. Of course, we could take the class of functionals $\mathcal{T}$ to be the set of evaluation functions $f(x)$ and so our approach does produce confidence bands too. It is easy to see, however, that without additional assumptions on the functions, these bands are hopelessly wide. We should also mention that another approach is to construct Bayesian posterior intervals as in Barber, Nason and Silverman (2002), for example. However, the frequentist coverage of such sets is unknown.

In Section 5 we gave a flavor of how information can be extracted from the confidence ball $\mathcal{C}_n$ using functionals. Beran (2000) discusses a different approach to exploring $\mathcal{C}_n$ which he calls "probing the confidence set." This involves plotting smooth and wiggly representatives from $\mathcal{C}_n$. A generalization of these ideas is to use families of what we call *parametric probes*. These are parameterized functionals tailored to look for specific features of the function such as jumps and bumps. In a future paper we will report on probes, as well as other practical issues that arise. In particular, we will report on confidence sets for other shrinkage schemes besides thresholding and linear modulators.

**Acknowledgments.** The authors thank the referees for helpful comments.

Department of Statistics
Carnegie Mellon University
Pittsburgh, Pennsylvania 15213
USA
e-mail: genovese@stat.cmu.edu
e-mail: larry@stat.cmu.edu